\RequirePackage{fix-cm}
\documentclass[smallextended]{svjour3}       
\smartqed  
\usepackage{graphicx}
\journalname{ }
\begin{document}

\title{Seven, eight, and nine mutually touching infinitely long straight round cylinders: Entanglement in Euclidean space
}

\titlerunning{Mutually touching cylinders}        

\author{Peter V. Pikhitsa         \and
       Mansoo Choi  
}


\institute{Peter V. Pikhitsa \at
              Seoul National University, 151-744 Seoul, Korea \\
              Tel.: +822-88-01694\\
              Fax: +822-88-06671\\
              \email{peter@snu.ac.kr}           
           \and
           Mansoo Choi \at
              Seoul National University,151-744 Seoul, Korea
}

\date{Received: date / Accepted: date}

\maketitle

\begin{abstract}
It has been a challenge to make seven straight round cylinders mutually touch before our now 10-year old discovery [Phys. Rev. Lett. {\bf 93}, 015505 (2004)] of configurations of seven mutually touching infinitely long round cylinders (then coined 7-knots). Because of the current interest in string-like objects and entanglement which occur in many fields of Physics it is useful to find a simple way to treat ensembles of straight infinite cylinders. Here we propose a treatment with a chirality matrix. By comparing 7-knot with variable radii with the one where all cylinders are of equal radii (here 7*-knot, which for the first time appeared in [phys. stat. solidi, b {\bf 246}, 2098 (2009)]), we show that the reduction of 7-knot with a set of non-equal cylinder radii to 7*-knot of equal radii is possible only for one topologically unique configuration, all other 7-knots being of different topology. We found novel configurations for mutually touching infinitely long round cylinders when their numbers are eight and ultimately nine (here coined 8-knots and 9-knots). Unlike the case of 7-knot, where one angular parameter (for a given set of fixed radii) may change by sweeping a scissor angle between two chosen cylinders, in case of 8- and 9-knots their degrees of freedom are completely exhausted by mutual touching so that their configurations are "frozen" for each given set of radii. For 8-knot the radii of any six cylinders may be changeable (for example, all taken equal) while two remaining are uniquely determined by the others. We show that 9-knot makes the ultimate configuration where only three cylinders can have changeable radii and the remaining six are determined by the three. Possible generalizations and connection with Physics are mentioned.

\end{abstract}

\section{Introduction}
\label{intro}
For centuries people combined different solid objects to construct other objects. On this way the restriction of a solid body being impenetrable by another body produces various packing problems. One example is packing spheres with a famous application to the liquid of hard spheres with the topological disorder still not understood. Intuitively, solid infinite cylinders should be more cumbersome because they can be entangled. Indeed, unlike compact bodies such as spheres, the infinite cylinders can touch one another at any distance, producing infinitely long spatial correlations. This may explain why for the whole history the first attempt to deal with the cylinders (finite rigid rods) was made by Lars Onsager only in 1949. It has long been known how to make mutually touch seven half-infinite cylinders, yet, the configuration essentially makes use of the possibility of three half-infinite cylinders to lie in one plane while touching both each other and the pivotal cylinder with their stumps. Placing other three cylinders right over them produces the configuration needed. However, being half-infinite is crucial for topology: a half-infinite cylinder does not "puncture" the space.  In 1986, in connection with the continuum percolation problem, we analyzed the possibility of arranging mutually touching arbitrary infinite cylinders invoking the method of degrees of freedom \cite{Ref1}. The method prompted only possible solutions to be found either numerically or analytically.  We predicted that seven mutually touching cylinders might exist and the configuration might not be rigid.

\begin{figure}[h]
  \includegraphics [width=0.85\textwidth]{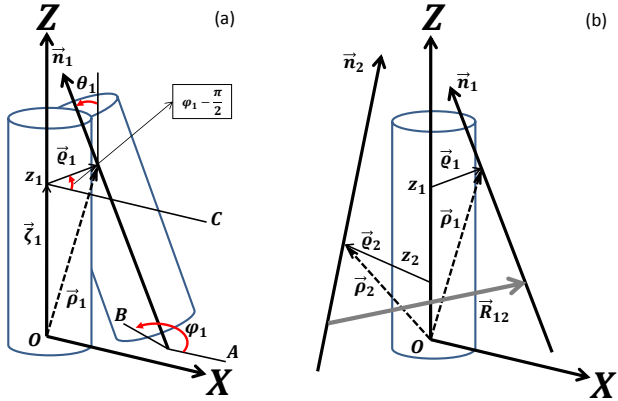}
\caption{ Geometry of contacts. The pivotal 0th cylinder lies along $Z$. $Y$ axis that extends to the right is not shown. (a) The 1st cylinder touching the pivotal vertical cylinder of unit radius $r_0=1$. Letters $A$ and $C$ mark the directions along $X$. (b) Two cylinders touching each other and the pivotal 0th cylinder }
\label{fig:1}       
\end{figure}

In 2004 we discovered \cite{Ref2} that the configuration of seven mutually touching straight round infinite cylinders indeed had a numerical solution (while solving 10 equations for 10 variables), which we named 7-knot. The configuration keeps its topology when one remaining degree of freedom (predicted in  \cite{Ref1} being the scissor angle $\theta_1$ (see Fig. \ref{fig:1}) between two chosen cylinders: the pivotal 0th cylinder (red) and the 1st cylinder (brown) in Fig. \ref{fig:2}), changes.  The rigorous nonlinear equations for the angles and positions (see below), which provide orientations of each cylinder to be touching all the others were solved with a renown Mathcad program with a computational accuracy of order of $10^{-13}$ (all relevant files with the texts of the  Mathcad programs as well as animations are given in  Online Resources 1-8)). The method was extended to regular lattices of mutually touching cylinders to produce an artificial matter with the ultimate Poisson ratio -1 because it could either shrink or expand as a whole due to cylinder entanglement \cite{Ref2}.

\begin{figure}[h]
  \includegraphics [width=0.85\textwidth]{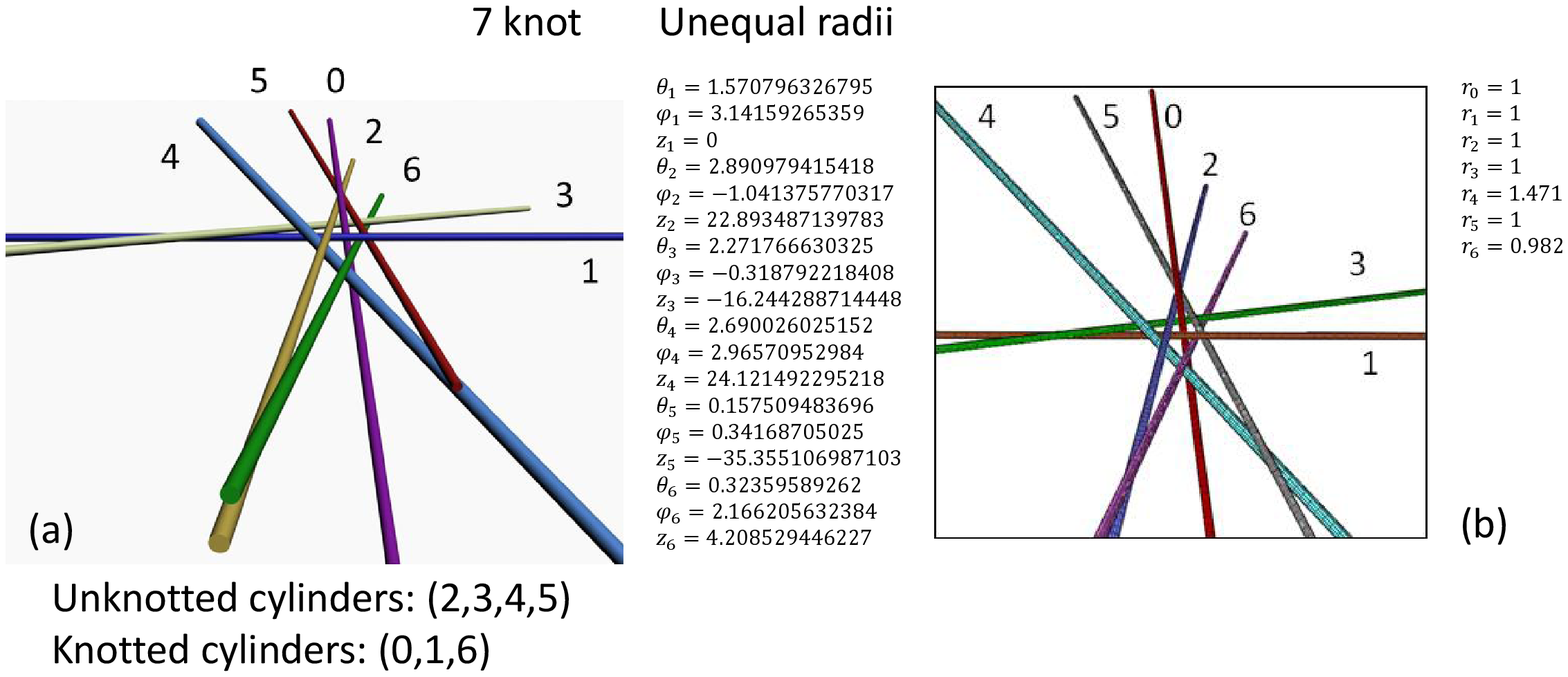}
\caption{7-knot of unequal cylinder radii. (a) The image of the configuration similar to the first time ever found configuration published in Fig. 1A of \cite{Ref2}; (b) the same configuration in perspective representation with the calculated parameters (the tolerance to zero was set $10^{-8}$). $0$ marks the pivotal cylinder and $1$ marks the first cylinder as given in Fig. \ref{fig:1}. Here the scissor angle was $\theta_1=\pi/2$. The configuration has $4$ unknotted cylinders each of which can be removed without disturbing other cylinders}
\label{fig:2}       
\end{figure}

\begin{figure}[h]
  \includegraphics [width=0.85\textwidth]{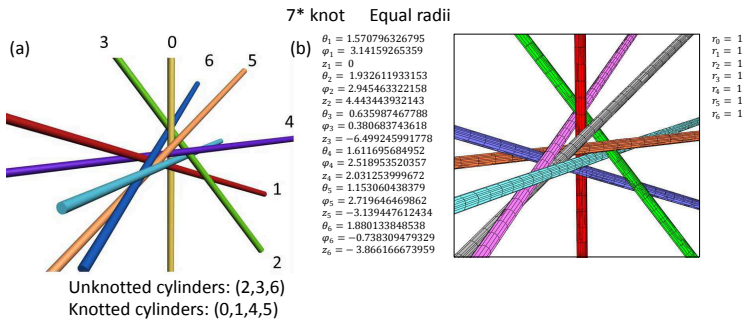}
\caption{Equal radii cylinders numbered. (a) The configuration similar to the one published for the first time ever in Fig. 1b of \cite{Ref3} (here the scissor angle was $\theta_1=\pi/2$); (b) the same as in (a) with calculated parameters (the tolerance to zero was set $10^{-8}$). The configuration has $3$ unknotted cylinders each of which can be removed to infinity without disturbing other cylinders. Therefore the configuration is topologically different from the configuration in Fig. \ref{fig:2}. The animation file with the sweeping angle $\theta_1$ is included in Online Resource 1 as well as the animation of the solid rotation of the configuration in Online Resource 6}
\label{fig:3}       
\end{figure}

In \cite{Ref2} the 7-knot was first found in a configuration topology that precluded having all the cylinders to be equal in radii (congruent). In 2009 we succeeded \cite{Ref3}in discovering other topologies and published the result with the 7-knot of the mutually touching infinite cylinders of equal radii (the congruent one, here equipped with asterisk: 7*-knot) (see Fig. \ref{fig:1}b) in \cite{Ref3} and Fig. \ref{fig:3}). However, probably because the results of constructing pure mathematical bodies were embedded into the physical context, the findings remained unknown for the mathematical community. A recent paper \cite{Ref4} rediscovered our finding of mutually touching seven infinite equal radii cylinders, the authors being unaware of our earlier results with cylinders of equal and non-equal radii and of the nuances of their topology.  It may illustrate a still existing gap between the worlds of Physics and Mathematics though both are intimately connected through the Plato's World \cite{Ref5}.

Meantime, the obvious difference between 7-knot and 7*-knot pushed us forward in curiosity questions about the topology of entanglement of the cylinders that constitute 7-knots, the limits of angles and radii where they exist, calculating the latitudinal and longitudinal angles while the scissor angle $\theta_1$ sweeps from the zero value to $\pi$ etc. Also there is the question of quantification of the difference in the topology between the configurations of the 7-knot with a set of non-equal radii that can be continuously reduced to 7*-knot with equal cylinders, and the basins of sets of non-equal radii that cannot be continuously reduced to the equal radii. The switching between the topologies always happens through infinity while the cylinders try to become parallel which is a kind of degeneracy. Therefore some orientational entanglement takes place which might be important to be understood on the basic level of straight infinite cylinders.

Here, we reveal the detailed structure of 7-knot along with our finding that natural extension to 8-and 9-knots does exist, as predicted in \cite{Ref2}. We computed the 8-knot and the ultimate 9-knot. These novel configurations of mutually touching round infinite cylinders lack the angular degrees of freedom because of having too many mutual contacts. It is not at all obvious that such rigid configurations should exist and simple comparison with the circles and spheres below illustrate the whole problem of the degrees of freedom for 7-, 8-, and 9- knots. Finally, we mention some possible consequences of cylinder entanglement for the growing field of various string-like objects in Physics and possible connection with Rodger Penrose's twistors \cite{Ref5}.

\section{Degrees of freedom}
\label{sec:1}
Let us begin with circles to explain the method of controlling configuration by the number $N_F$ of degrees of freedom. For N circles on the plane being in mutual contacts the number $N_F$ satisfies the degree of freedom equation (DFE):
\begin{equation}
N_F  = 2N-N(N-1)/2-3.
\label{eq:1}
\end{equation}
Indeed, each circle has $2$ degrees of freedom ($x$,$y$ positions of its center on the plane), that is $2N$ in total for all the freedom and $-N(N-1)/2$ stands for pairwise restrictions due to mutual contacts. The last term is the number of solid degrees of any two dimensional configuration as a whole: it can be rotated $(-1)$ and translated $(-2)$. Easy to see that for $N=3$ Eq. (\ref{eq:1}) gives $N_F=0$. This rigid configuration actually exists for any radii: one can always mutually contact $3$ circles. If we restrict ourselves with outside contacts then we would have a single configuration with a mirror counterpartner. Anyway mirror degrees do not count here because they are discrete. Then let us take $N=4$ which gives $N_F=-1$. Having $N_F=-1$ indicates that the fourth circle should have a predetermined radius; others are changeable (so-called Apollonius circles). The next case is more dramatic: $N=5$ gives us $N_F=-3$. This should be enough because we still have $5$ radii to be adjusted as we please (more precisely, four radii, because one radius should be a measuring unit). However, such an ultimate (one cannot put $N=6$ because $N_F=-6$ is out of range, only $5$ radii out of $6$ can be changeable) configuration does not exist. It means that even if a configuration is not forbidden by the counting of the degrees of freedom, it is not at all obvious that it can exist.

For spheres instead of Eq. (\ref{eq:1}) one has
\begin{equation}
N_F= 3N-N(N-1)/2-6.
\label{eq:2}
\end{equation}
For $N=4$,$N_F=0$, which is also rigid at some changeable radii. For N=5 Eq. (\ref{eq:2}) gives $N_F=-1$, which is quite realizable provided that one sphere radius is determined by the rest four. However, $N=6$ gives $N_F=-3$ that does not exist though we still have enough radii to fix. The ultimate limit would be $N=7$ with $N_F=-6$, which is still enough though also does not exist.

Let us see what occurs for infinite cylinders. Here the situation is different: all configurations predicted by DFE exist. We posed this question with the connection to the topology of the percolation cluster in \cite{Ref1} where the corresponding DFE was discussed (see also \cite{Ref2} with more relevant examples). The equation for N cylinders is \cite{Ref1}, \cite{Ref2}:
\begin{equation}
N_F  = 4N-N(N-1)/2-6.
\label{eq:3}
\end{equation}
Eq. (\ref{eq:3}) looks like a minor alternation as compared with Eq. (\ref{eq:2}) but it takes into account the crucially important affine degree of freedom, that for infinite cylinders implies that a cylinder axis can always punch a two dimensional plane at a point ($2$ degrees for $x$,$y$) and can be ray-like oriented ($2$ degrees for two angles $\theta_1$ and $\varphi_1$, see Fig. \ref{fig:1}). Eq. (\ref{eq:3}) always predicts actually existing configurations. One of them was 7-knot \cite{Ref2}, \cite{Ref3},\cite{Ref4} and the other two will be shown below to be 8-knot and 9-knot. Let us use Eq. (\ref{eq:3}) for $N=7$. From Eq. (\ref{eq:3}) we have $N_F=1$ which means that 7-knot is not rigid: one degree (it can be the scissor angle $\theta_1$) remains. The calculations indeed proved this \cite{Ref2}, \cite{Ref3},\cite{Ref4}. For $N=8$ we obtain $N_F=-2$ so two radii of the cylinders should be determined and they are, as we show below. Then, $N=9$ gives $N_F=-6$ with only $3$ cylinders of changeable radii (again only two in reality as far as one radius always goes for a measuring unit). This ultimate configuration exists, in spite of the non-existence of ultimate configurations for circles and spheres.

We write down the equations for the mutual contacts in the vector notations. The number of equations we choose to solve is $(N-1)(N-2)/2$ (see below) which makes $15$, $21$, and $28$ in case of 7-, 8-, and 9-knot, respectively.  We solve them numerically with the Mathcad solver (all relevant texts of the programs given in Online Resources 2-5). The sets of equations also could be simplified more as to their number and variables. For example, in \cite{Ref2} after some analytical work we used only $10$ equations to fix topology and to solve only for angles but in \cite{Ref3} we used either $10$ or $15$ equations with the same result (see below). Our practice shows that the solution converges nearly as fast even for 9-knot with $28$ equations, or for 8-knot with $21$ equations, as for 7-knot with either $15$ or $10$ equations depending on how one defines the variables.

\section{Basic equations}
\label{sec:2}

The basic equations for the contacting cylinders can be obtained from Fig. \ref{fig:1}. We use the 0th cylinder of the radius $r_0=1$, taken to be a measuring unit, as the pivotal one which axis lies along $Z$. All other cylinders are constructed as touching it by default. The 1st cylinder of radius $r_1$ touching the 0th cylinder is described by the unit vector $\vec n_1$ that defines the orientation of the cylinder and the vector $\vec \varrho_1$ connecting the origin with a point on the cylinder axis. The longitudinal angle $\varphi_1$ is counted counterclockwise from the direction of $X$ axis, marked by letter $A$, to the direction of the vertical projection of the cylinder axis onto the plane $XY$, marked by letter $B$ in Fig.\ref{fig:1}a. The axis $Y$, not shown in Fig.  \ref{fig:1}, extends to the right. The latitudinal angle $\theta_1$ is shown in Fig.  \ref{fig:1}a. Then we define

\begin{equation}
\vec n_1 = \left( \begin{array}{c}
 \sin(\theta_1 )\cos(\varphi_1) \\
\sin(\theta_1 )\sin(\varphi_1) \\
  \cos(\theta_1 )
\end{array} \right),
\label{eq:4}
\end{equation}

\begin{equation}
\vec\rho_1 = \left( \begin{array}{c}
(r_1+r_0 )\sin(\varphi_1 ) \\
-(r_1+r_0 )  \cos(\varphi_1  ) \\
 z_1
\end{array}\right).
\label{eq:5}
\end{equation}

In such a representation it is easy to check from Eqs. (\ref{eq:4}), (\ref{eq:5}) that the vector
\begin{equation}
\vec \varrho_1 = \vec \rho_1 - \vec \zeta_1=  \left( \begin{array}{c}
(r_1+r_0 )\sin(\varphi_1 ) \\
-(r_1+r_0 )  \cos(\varphi_1  ) \\
 z_1 
\end{array} \right)
-  \left( \begin{array}{c}
0 \\
0 \\
 z_1 
\end{array}\right)
=   \left( \begin{array}{c}
(r_1+r_0 )\sin(\varphi_1 ) \\
-(r_1+r_0 )  \cos(\varphi_1  ) \\
0 
\end{array}\right)
\label{eq:6}
\end{equation}
which length is $(r_1+r_0 )$ (Fig. \ref{fig:1}a) is orthogonal to both $Z$ axis and $\vec n_1$, therefore the 1st cylinder is indeed touching the 0th cylinder by construction. Any point on the axis of the 1st cylinder is described by the vector $\vec \rho_1+\vec n_1 T_1$, where $T_1$ is the affine parameter along the cylinder axis. Thus, as far as we fixed the pivotal cylinder to be vertical, all other cylinders are now characterized by only 3 degrees of freedom each (provided we keep the radii fixed): ($\theta_1$,$\varphi_1$,$z_1$ ),($\theta_2$,$\varphi_2$,$z_2$  ), etc. To fix the configuration completely, we always keep $z_1=0$ and $\varphi_1=\pi$, which removes 2 degrees of freedom. Now we have $(N-1)$ cylinders to make touch mutually. That allows one to rewrite DFE of Eq. (\ref{eq:3}) in the form:
\begin{equation}
N_F  = 3(N-1)-(N-1)(N-2)/2-2,
\label{eq:7}
\end{equation}
where $N_F$ is of course identical to the one from Eq. (\ref{eq:3}). Yet instead of $21$ contact equations for $N=7$ we are having $15$. For $N=8$ and $9$ they are $21$ and $28$, respectively.

Let us turn to Fig. \ref{fig:1}b and write down the equations of contact between the cylinders. Vector $\vec R _{12}$ connecting a point on one cylinder axis with a point of the other cylinder axis, can be written as:
\begin{equation}
\vec R_{12}= \vec \rho_1+\vec n_1 T_1- (\vec \rho_2+\vec n_2 T_2 ),
\label{eq:8}
\end{equation}
where again $T_1$ and $T_2$ are the affine parameters along the cylinder axes. As far as at the touching point $\vec R_{12}$  being the shortest vector should be orthogonal to both axes of the cylinders $1$ and $2$, one can write down two equations to pin-point the place of touching which determine $T_1$ and $T_2$:
\begin{equation}
 \begin{array}{c}
\vec R_{12} \vec n_1=0\\
\vec R_{12} \vec n_2=0 \end{array}  \\ .
\label{eq:9}
\end{equation}

Then for $T_1$ and $T_2$ one obtains from Eq. (\ref{eq:9}):
\begin{equation}
T_1=\frac{-(\vec n_1 (\vec \rho_1- \vec \rho_2 ) )+(\vec n_1 \vec n_2 )(\vec n_2 (\vec \rho_1- \vec \rho_2  ) )}{1-(\vec n_1 \vec n_2 )^2 }  ,
\label{eq:10}
\end{equation}

\begin{equation}
T_2=\frac{-(\vec n_1 \vec n_2 )(\vec n_1 (\vec \rho_1- \vec \rho_2 ) )+(\vec n_2 (\vec \rho_1- \vec \rho_2  ) )}{1-(\vec n_1 \vec n_2 )^2 }  .
\label{eq:11}
\end{equation}
On substituting Eqs. (\ref{eq:10}) and (\ref{eq:11}) into Eq. (\ref{eq:8}) one finally gets the first equation to solve:

\begin{equation}
(\vec R_{12}  )^2=(r_1+r_2 )^2 
\label{eq:12}
\end{equation}
with all the rest equations being obtained by replacing $(\theta_1 \rightarrow \theta_i,\theta_2 \rightarrow \theta_j )$,$(\varphi_1\rightarrow\varphi_i,\varphi_2\rightarrow\phi_j )$,$(z_1\rightarrow z_i,z_2\rightarrow z_j)$,$(r_1\rightarrow r_i,r_2\rightarrow r_j)$, where $i=1\dots (N-1)$,$j=1\dots (N-1)$ and $i<j$.  The solving block that in Mathcad begins with "Given" and ends with "Find" contains either $15$, $21$, or $28$ equations for 7-, 8-, and 9-knots, correspondingly, in accordance with Eq. (\ref{eq:7}). The accuracy of the Mathcad calculations was $10^{-13}$  in order. The equations were calculated for angles and sometimes for radii when it was appropriate.
In case of $N=7$ the numbers of unknowns was $3(N-1)-3=15$, thus coinciding with the number of the equations, because we fixed $\theta_1$ (any angle from $0$ to $\pi$), $\varphi_1=\pi$, and $z_1=0$. In case of $N=8$ the number of unknowns was $3(N-1)-2+2=21$ because we now calculated also for $\theta_1$, remaining $\varphi_1=\pi$, and $z_1=0$ fixed, and also two of the radii were taken as unknowns and calculated, for example $r_4$ and $r_7$ as in Fig. \ref{fig:8}a below, or $r_5$ and $r_7$  in Fig. \ref{fig:8}b. In case of $N=9$ the number of unknowns was $3(N-1)-2+6=28$, where we added $6$ radii $r_3$,$r_4$,$r_5$,$r_6$,$r_7$,$r_8$ to unknowns.

We could additionally control the calculations and topology of configurations by utilizing the possibility to calculate $z_i$ also analytically as it was performed in [2], [3]. Let us fix $z_1=0$. It is easy to make sure from Eqs. (\ref{eq:5}), (\ref{eq:10}), (\ref{eq:11}) that $\vec R_{12}$ of Eq. (\ref{eq:8}) is a linear function of $(z_2-z_1 )=z_2$. After rewriting Eq. (\ref{eq:8}) as

\begin{equation}
\vec R_{12} (z_2 )=\vec \Delta_{12}+z_2 \vec \Lambda_{12},  
\label{eq:13}
\end{equation}
from Eq. (\ref{eq:13}) one obtains by inserting $z_2$ either $0$ or $1$: $\vec \Delta_{12}=\vec R_{12} (0)$ and $\vec \Lambda_{12}=\vec R_{12} (1)-\vec R_{12}(0)$.  Substituting Eq. (\ref{eq:13}) into Eq. (\ref{eq:12}) one obtains a quadratic equation for $z_2$
\begin{equation}
(\vec \Lambda_{12} )^2 z_2^2  + 2(\vec \Lambda_{12} \vec \Delta_{12}  ) z_2+(\vec \Delta_{12} )^2=(r_1+r_2)^2,
\label{eq:14}
\end{equation}
which solution is:
\begin{equation}
z_2=\frac{-(\vec \Lambda_{12} \vec \Delta_{12} )+\eta_2 \sqrt{(\vec \Lambda_{12} \vec \Delta_{12}  )^2+[(r_2+r_1 )^2-(\vec \Delta_{12} )^2 ] (\vec \Lambda_{12} )^2 }}{(\vec \Lambda_{12} )^2}   .
\label{eq:15}
\end{equation}
Here $\eta_2$ defines the sign of the square root to be either $\eta_2=1$ or $\eta_2=-1$ and discriminates the topology. Replacing index $2\rightarrow i$ gives the solutions for $i=2,3,\dots $.The choice of topology for 7-knot in Fig. \ref{fig:2} corresponds to $\eta_2=-1$; $\eta_3=-1$; $\eta_4=1$; $\eta_5=1$; $\eta_6=1$. For 7*-knot of Fig. \ref{fig:3} the topology corresponds to $\eta_2=1$; $\eta_3=1$; $\eta_4=1$; $\eta_5=-1$; $\eta_6=-1$, the same set as in \cite{Ref3}. With the help of Eq. (\ref{eq:15}) we can compare the precision of the calculations. For example, for 7-knot in Fig. \ref{fig:2} we obtained (setting the tolerance to zero at $10^{-8}$) $z_i$ in two ways: once calculating with full $15$ equations similar to Eq. (\ref{eq:12}) for all $15$ $\theta_i$,$\varphi_i$,$z_i$  (the second column in Table \ref{tab:1}) and then re-calculating $z_i$ analytically from Eq. (\ref{eq:15}) (the third column in Table \ref{tab:1}). One can see rather good accuracy for both 7-knot and 7*-knot.
\begin{table}[h]
\caption{Comparison between computationally and analytically obtained parameters $z_i$.}
\label{tab:1}       
\begin{tabular}{llllll}
\hline\noalign{\smallskip}
7-knot & Calc. with Eq. (12) & Calc. with Eq. (15)& 7*-knot & Calc. with Eq. (12) & Calc. with Eq. (15) \\
\noalign{\smallskip}\hline\noalign{\smallskip}
$z_2$ & -22.893487139783 & -22.893487139783 &	$z_2$ &	4.443443932143	& 4.443443932142 \\
$z_3$ &	-16.244288714448 &	-16.244288714448 &	$z_3$ &	-6.499245991778 &	-6.499245991776\\
$z_4$ &	24.121492295218	 & 24.121492295218 &	$z_4$ &	2.031253999672 &	2.031253999672\\
$z_5$ &	-35.355106987103 &	-35.355106987102 &	$z_5$ &	-3.139447612434 &	-3.139447612434\\
$z_6$ &	4.208529446227	 & 4.208529446227 &	$z_6$	& -3.866166673959 &	-3.866166673959\\
\noalign{\smallskip}\hline
\end{tabular}
\end{table}

The main challenge in the calculations for $N=7$ was first to locate the appropriate initial values for $\theta_i$,$\varphi_i$,$z_i$,$r_i$ to find any solution possible which was realized in \cite{Ref2} (Fig. \ref{eq:2}). The second challenge was to find the topology that provides the most entangled structure with only three "unknotted" cylinders (the cylinders any of which  can be translated to infinity without disturbing all the rest, see below) like the three cylinders, marked $6$(violet), $3$(green), $2$(blue) in 7*-knot with all cylinders equal, in Fig. \ref{fig:3}b.  Understanding that having the tightest configuration might be important for looking for configurations with $N>7$ stipulated our success in finding firstly 8-knot and then the ultimate 9-knot. Let us describe the results obtained.

\section{The 7-knot}
\label{sec:3}
As it is seen from Eqs. ( \ref{eq:3}) and (\ref{eq:4}), $N_F=1$ for $N=7$. This degree of freedom is controlled by the angle $\theta_1$ which can be sweeping from $0$ to $\pi$ while the topology of a configuration does not change.  However, the topology is connected with the set of radii. For example, one can have the 7-knot from Fig. \ref{fig:2} with the radii $(r_1,r_2,r_3,r_4,r_5,r_6)=(1,1,1,1,2,1)$ as well as the ~7*-knot (tilde indicates the same topology but not all equivalent radii) from Fig. \ref{fig:3} with the same set $(1,1,1,1,2,1)$, so that $r_5=2$ for both cases. 
Let us vary only $r_5$. Then, while for ~7*-knot $r_5$ can be reduced to $1$, in case of 7-knot $r_5$ can be reduced only to $r_5\approx 1.265$, where some of the cylinders try to get parallel. On the other hand, ~7*-knot topology configuration cannot exist for $r_5>3.55$ as the calculations show (see Fig. \ref{fig:6}a for $r_5=3$ as the cylinders try to get parallel). A traverse in the $6D$ space of radii set thus switches the configuration to different topology. Topology is preserved in radii variations until the cylinders become parallel. There should exist non-overlapping sets of radii ("basins" or "valleys") that have one-to-one correspondence with the topology. At least for the vicinity of the set of the radii $(1,1,1,1,1,1)$ it is true: this set selects 7*-knot topology uniquely.

\begin{figure}[h]
  \includegraphics [width=0.85\textwidth]{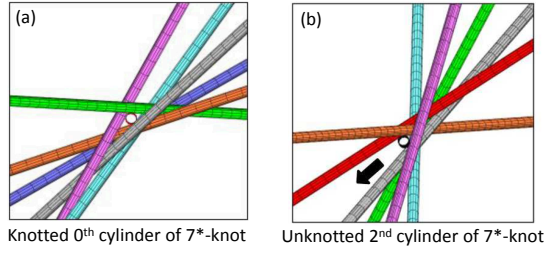}
\caption{Two different projections of 7*-knot illustrate the notion of knotted (a) and unknotted (b) cylinders. The black arrow shows how the 2nd cylinder can be removed to infinity without disturbing the rest of the cylinders }
\label{fig:4}       
\end{figure}

Let us show quite simply that the topology of the two configurations in Fig. \ref{fig:2} and Fig. \ref{fig:3} is different with the help of the notions of "knotted" and "unknotted" cylinders, illustrated in Fig. \ref{fig:4} on example of two projections of 7*-knot. Indeed, it is impossible to remove the 0th cylinder by translating it in any direction orthogonal to its axis without disturbing the rest in Fig. \ref{fig:4}a. On contrary, the 2nd cylinder can be translated to infinity in the direction of the black arrow in Fig. \ref{fig:4}b. One can see by inspection for the 7-knot in Fig. \ref{fig:2} that there are four unknotted cylinders $2$, $3$, $4$, $5$ and for 7*-knot in Fig.  \ref{fig:3} there are only three unknotted cylinders $2$, $3$, $6$; therefore the configurations are topologically different.

Note that even if we "relax" 7-knot configurations by translating the cylinders arbitrarily in directions orthogonal to their axes (so that they stop touching each other), but without changing their orientations, the configurations would retain their entanglement and topology as far as the topology can change only through the parallelism of cylinders. It indicates that to any configuration of infinite cylinders some numerical characteristics can be ascribed which is in one to one correspondence with the topology of the configuration. We suggest the simplest symmetric matrix that discerns the topology with a typical normalized element of chirality:
 \begin{equation}
P_{12}=\frac{(\vec n_1\times \vec n_2 ) \vec R_{12}}{\left|(\vec n_1\times \vec n_2 ) \vec R_{12} \right|} ,  \label{eq:16}
\end{equation}
which renders either $+1$ or $-1$, and we put zero on the diagonal. For the orientation of vectors $\vec n_1$ and $\vec n_2$ being defined with respect to the pivotal cylinder as shown in Fig. \ref{fig:1}, and vector $\vec R _{12}$ being orthogonal to both direction while satisfying Eq. (\ref{eq:8}), any change from $+1$ to $-1$ may occur either through the parallel configuration (which is a degeneracy) or through the mirror reflection (which is a discreet operation). Thus, chirality matrix $P_{ij}$ of Eq. (\ref{eq:16}) is preserved when the topology of the entanglement of the cylinders is preserved, and it makes the classification of entangled configurations of straight cylinders much simpler than classification of conventional one-component knots. For 7-knot from Fig. \ref{fig:2}, the matrix $P_{ij}$ is calculated to be
 \begin{equation}
P_{ij}=\left( \begin{array}{ccccccc}
0&+1&+1&+1&+1&+1&+1\\
+1&0&+1&+1&+1&+1&+1\\
+1&+1&0&-1&+1&+1&-1\\
+1&+1&-1&0&-1&+1&-1\\
+1&+1&+1&-1&0&-1&-1\\
+1&+1&+1&+1&-1&0&-1\\
+1&+1&-1&-1&-1&-1&0 \end{array} \right)
\label{eq:17}
\end{equation}
with its determinant $|P_{ij} |=10$. For 7*-knot given in Fig. \ref{fig:3}, the matrix $P_{ij}$ is obviously different
 \begin{equation}
P_{ij}=\left( \begin{array}{ccccccc}
0&+1&+1&+1&+1&+1&+1\\
+1&0&+1&+1&+1&-1&+1\\
+1&+1&0&-1&-1&-1&+1\\
+1&+1&-1&0&-1&+1&+1\\
+1&+1&-1&-1&0&+1&-1\\
+1&-1&-1&+1&+1&0&+1\\
+1&+1&+1&+1&-1&+1&0 
\end{array} \right)
\label{eq:18}
\end{equation}
 with its determinant $|P_{ij} |=-18$.  We made sure that the matrix $P_{ij}$ remains invariant while the scissor angle $\theta_1$ sweeps the whole range from $0$ to $\pi$. The first row and the first column always contain only $+1$ because of the chosen right-hand orientation way for angles $\theta_i$  and $\varphi_i$  (Fig. \ref{fig:1}a) at the contact with the pivotal 0th cylinder. Given this matrix one can unambiguously interweave 7 cylinders in a 7-knot. There are also interesting properties of the matrix of Eq. (\ref{eq:16}) with non-normalized elements $\mathcal{P}_{ij}=(\vec n_i\times \vec n_j ) \vec R_{ij}$.
 
7-knots may be connected with 7-tangles with some reservations. Indeed, according to Conway's definition, an n-tangle is a proper embedding of a disjoint collection of $n$ arcs into a 3-ball when the embedding must send the endpoints of the arcs to 2n marked points on the ball's boundary. For the case of 7-knot the endpoints are strictly polar points of the ball (Fig. 5a) where the position of $2n$ points is dictated by the radii and the conditions of touching. Yet the arcs, equipped with non-zero thickness while representing the touching cylinders, should also be pairwise cross-linked at one point of crossing. This cross-link makes it impossible for a third arc to slip "between" the two arcs, thus producing a "knotting" in the 7-tangle.

\begin{figure}[h]
  \includegraphics [width=0.85\textwidth]{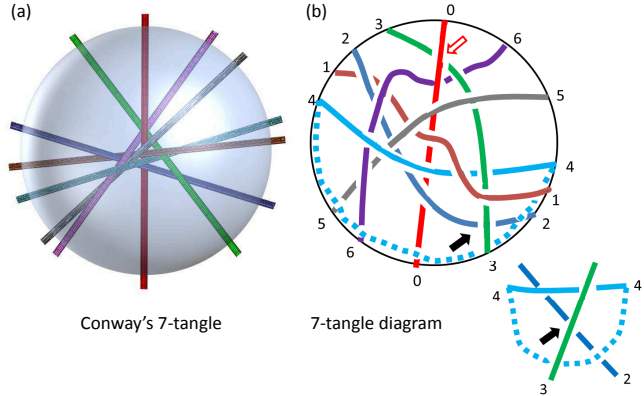}
\caption{Correspondence of 7*-knot with Conway's 7-tangle. (a) The cylinders punch the surface of the 3-ball ("celestial sphere") in opposite polar points; (b) the projection 7-tangle diagram of (a) equipped with contacts of arcs: the black arrows show one of the contacts that makes the 4th (cyan) arc knotted.  The 4th arc cannot be re-deformed into the dotted one that runs along the great circle. From the inset below it is understood why. The contact is impassible for the 4th cyan close contour that encompasses it neither for shrinking to null nor for putting off }
\label{fig:5}       
\end{figure}

\begin{figure}[h]
  \includegraphics [width=0.85\textwidth]{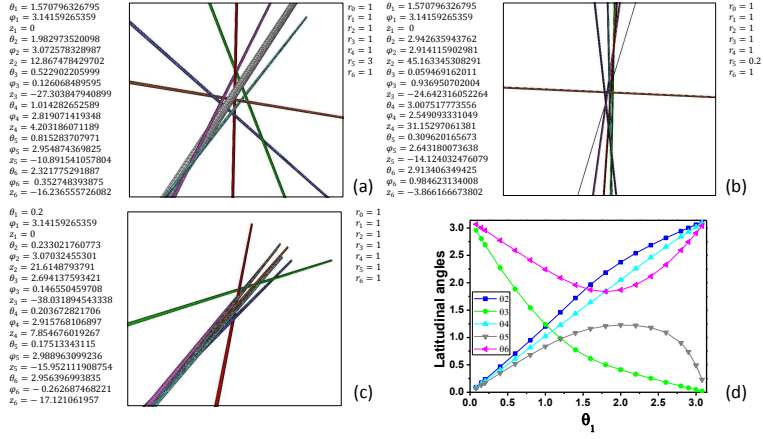}
\caption{Illustration of the evolution of so-called ~7*-knot of the same topology as 7*-knot of equal radii when parameters change. (a) As $r_5$ grows larger some cylinders tend to become parallel; (b) when $r_5$ grows smaller some cylinders also tend to become parallel; (c) for 7*-knot as $\theta_1\rightarrow 0$ cylinders grow parallel, the same holds for $\theta_1\rightarrow \pi$ (see Online Resource 1 animation); (d) the graph demonstrates the dependence of latitudinal angles $\theta_i$ on the scissor angle $\theta_1$ for 7*-knot. The curves are attracted either to $0$ or to $\pi$ as $\theta_1$ approaches either  $0$ or $\pi$ (the cylinders get parallel)}
\label{fig:6}       
\end{figure}

Let us illustrate it by a projection diagram in Fig. \ref{fig:5}b. It is known that without loss of generality, one can consider the marked points on the 3-ball boundary to lie on a great circle. The tangle can be arranged to be in general position with respect to the projection onto the flat disc bounded by the great circle which is the tangle diagram shown in Fig. \ref{fig:5}b. The inset to this figure illustrates in detail the attempt to re-deform the 4th cyan arc into the one shown with the dotted line which lies along the great circle and symbolizes the movement of the corresponding cylinder to infinity. It is impossible because at the position of the cross-link (indicated by the black arrow) between the 2nd (blue) and the 3rd (green) arcs, the 4th arc should "slip" between the 2nd and the 3rd arcs, which is prevented by their cross-link (see the caption to Fig. \ref{fig:5}b). On contrary, the 6th (violet) arc can be re-deformed in the North-West direction to lie along the great circle. The 6th arc, while passing through the cross-link marked by the red arrow, never gets between the 0th (red) and the 3rd (green) arcs.  Therefore, in consistency with above-said, the 4th arc (as well as the 4th cylinder) is "knotted" and the 6th arc (as well as the 6th cylinder) is "unknotted". The 7-tangle diagram with cross-links indeed reflects the topology of the contacting cylinders in the 7*-knot.

Let us give some computational results for 7*-knot (the tolerance to zero was set $10^{-6}$). As we mentioned above, the scissor angle $\theta_1$ between the 0th pivotal cylinder (red) and the 1st cylinder (brown) varies from $0$ to $\pi$. In Fig. \ref{fig:6}c we give a side view of the configuration where $\theta_1$ approaches $0$ and the cylinders get parallel. Fig. \ref{fig:6}d gives the computed dependence of all latitudinal angles $\theta_i$ vs $\theta_1$. One can notice that the 6th and the 5th cylinders do not upend from $0$ to $\pi$ or vise versa.

\begin{figure}[h]
  \includegraphics [width=0.85\textwidth]{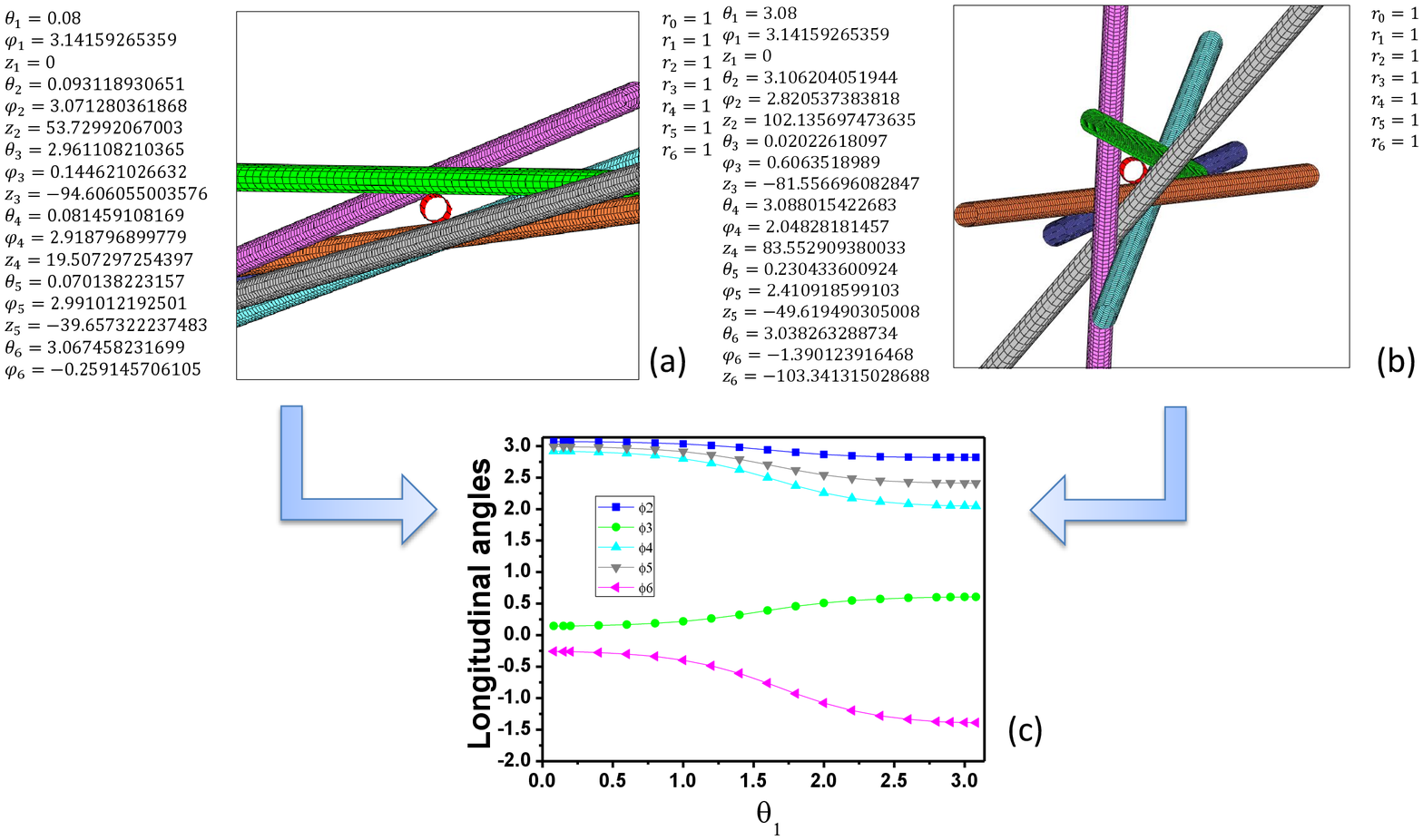}
\caption{Two limiting near parallel configurations as being viewed along the 0th pivotal cylinder with (a) $\theta_1\rightarrow 0$ and (b) $\theta_1\rightarrow \pi$; (c) The graph demonstrates the dependence of longitudinal angles $\varphi_1$ on the scissor angle $\theta_1$. The blue arrows indicate the limiting angles $\varphi_i$ for configurations in (a) and (b)    }
\label{fig:7}       
\end{figure}

 Being projected along the pivotal cylinder the 7*-knot demonstrates two limiting projections: one in Fig. \ref{fig:7}a corresponds to $\theta_1\rightarrow 0$ and the other in Fig. \ref{fig:7}b corresponds to $\theta_1\rightarrow \pi$. The total computed dependencies of $\varphi_i$ vs $\theta_1$ are given in Fig. \ref{fig:7}c.


\section{The 8-knot}
\label{sec:4}
\begin{figure}[h]
  \includegraphics [width=0.85\textwidth]{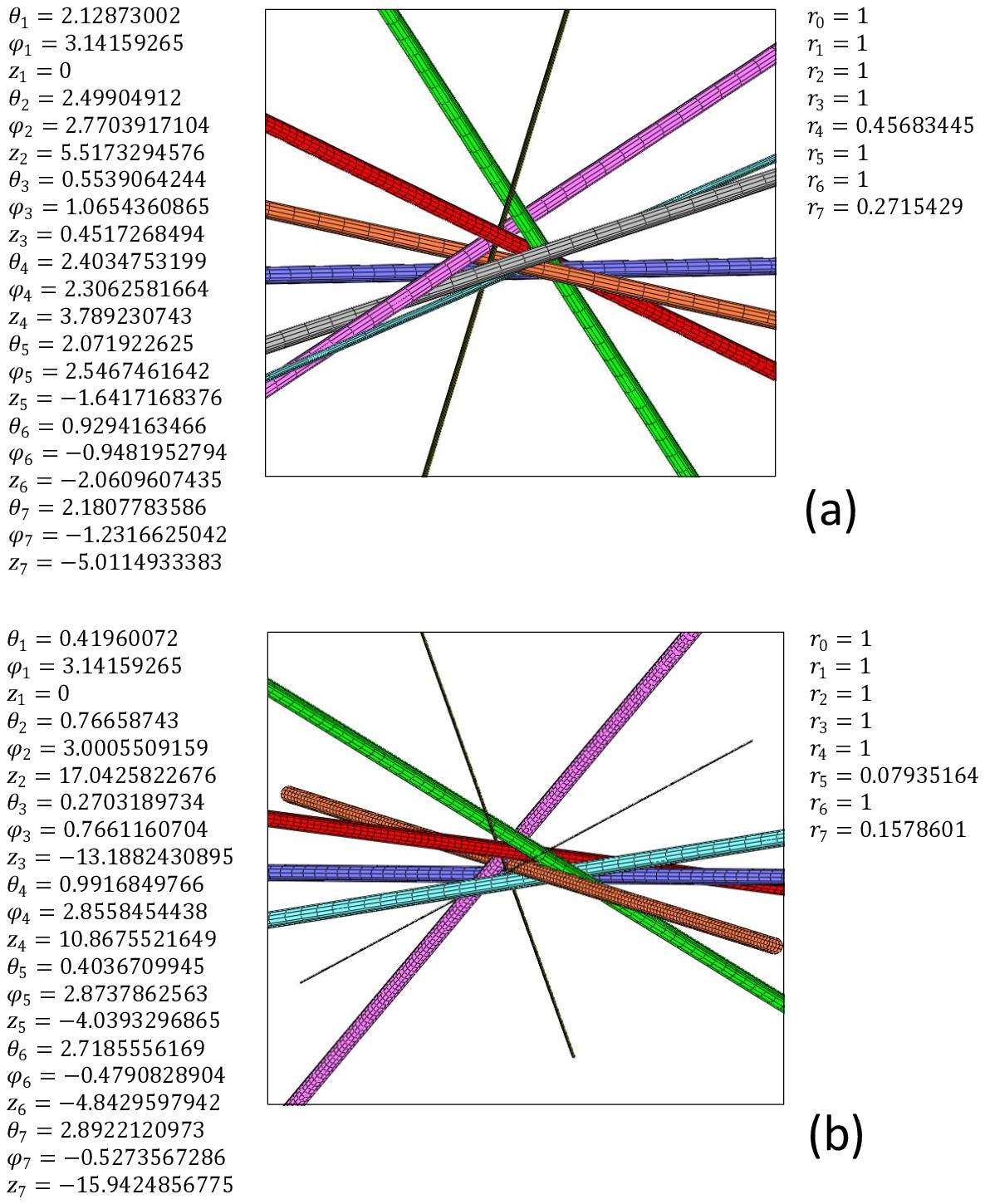}
\caption{Two configurations of 8-knot with 6 cylinders of unit radii; the remaining two radii are chosen to be calculated: (a) $r_4$ and $r_7$ are chosen; (b) $r_5$ and $r_7$ are chosen   }
\label{fig:8}       
\end{figure}

Although it is possible to calculate 8-knot with changeable nonequivalent radii we present in Fig. \ref{fig:8} only two unique configurations with calculated ($r_4=0.45683445$; $r_7=0.2715429$) and ($r_5=0.07935164$; $r_7=0.1578601$) while all the rest $6$ cylinders in the 8-knot are of the unit radii.
Any of such configurations with two determined radii and all the rest of the unit size is unique, being $28$ configurations in total if one combines all possible pairs of radii. They can be calculated by using the Mathcad program text given in Online Resource 4.  The 8-tangle for this configuration can be drawn as well. The control of the topology of the configuration of the cylinders according to their degree of knotting is the same as for 7-knot. For example, exactly like in the 7-knot with equal cylinders only three cylinders 6th (violet), 3rd (green), and 2nd (blue) are unknotted. Finally, the matrix $P_{ij}$  for both configurations in Fig. \ref{fig:8}  is
\begin{equation}
P_{ij}^{(8)}=\left( \begin{array}{cccccccc}
0&+1&+1&+1&+1&+1&+1&+1\\
+1&0&+1&+1&+1&-1&+1&+1\\
+1&+1&0&-1&-1&-1&+1&+1\\
+1&+1&-1&0&-1&+1&+1&-1\\
+1&+1&-1&-1&0&+1&-1&-1\\
+1&-1&-1&+1&+1&0&+1&-1\\
+1&+1&+1&+1&-1&+1&0&-1\\
+1&+1&+1&-1&-1&-1&-1&0 
\end{array} \right)
\label{eq:19}   
\end{equation}
with its determinant $|P_ij^{(8)} |=9$. We introduced the upper script $(8)$ to indicate the 8-knot. One can notice that the left upper corner $7\times 7$ block of the matrix exactly coincides with the matrix in Eq. \ref{eq:18}. Indeed, because we constructed the 8-knot by making use of the tightest 7*-knot configuration to insert the eighth cylinder, this topology remained "frozen".

\section{The 9-knot}
\label{sec:5}

\begin{figure}[h]
  \includegraphics [width=0.85\textwidth]{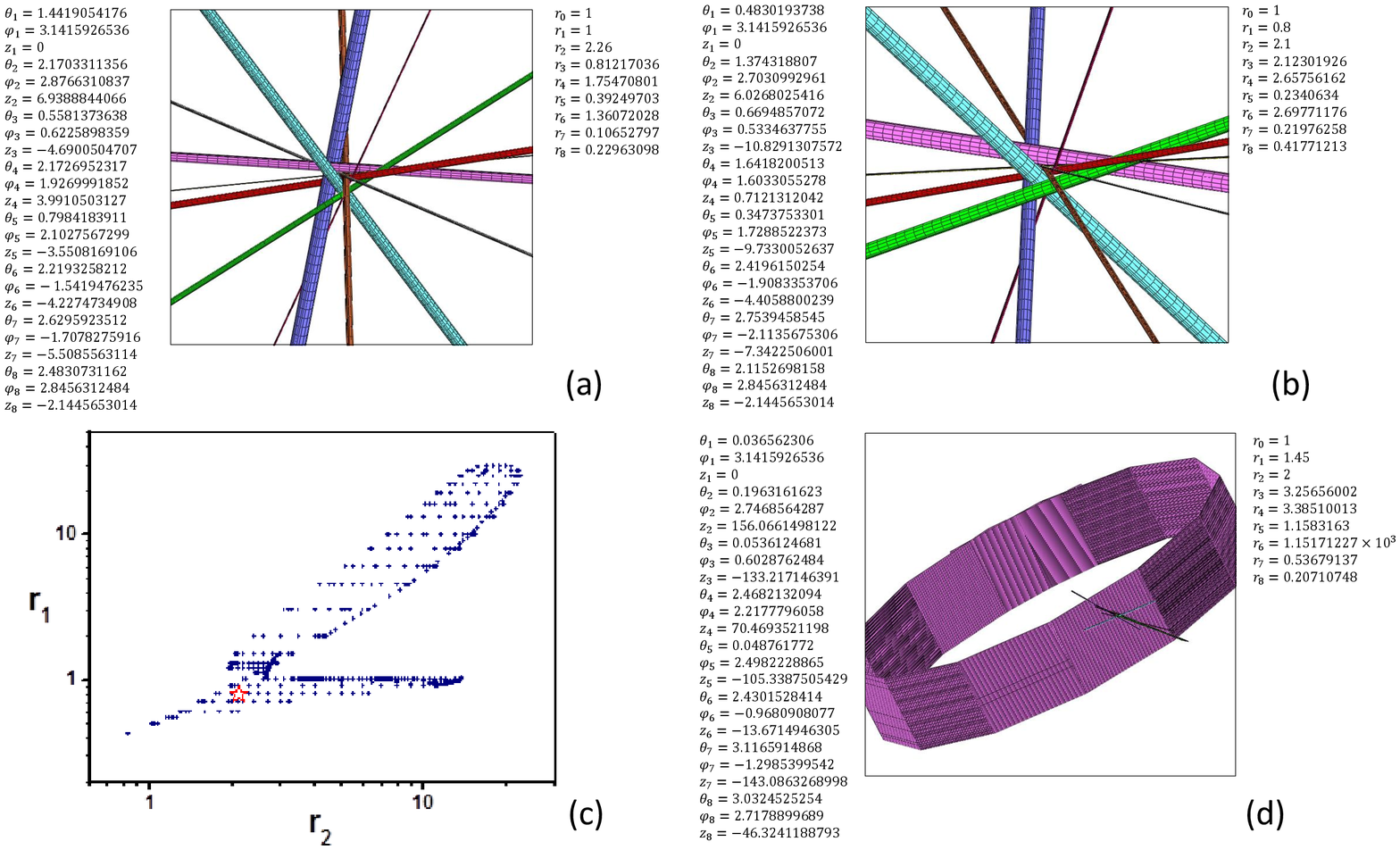}
\caption{Various configurations for 9-knot. (a) As obtained for the first time; (b) the one having the least ratio between the largest and the smallest radii; (c) the explored domain of parameters ($r_1$, $r_2$) for existing configurations in double-log scale. 587 dots mark existing configurations. The configuration (b) is marked with the star; (d) one of the configurations which is near the border of the domain showing the tendency of one radius to infinity (the part of the huge cylinder looking like a ring) with some others being nearly parallel    }
\label{fig:9}       
\end{figure}

This one is more interesting because it is the ultimate configuration possible for the infinite round cylinders. It closely resembles a complicated version of Apollonius circles in the respect that three radii determine all the rest. We present two of 9-knot configurations of the same topology in Fig. \ref{fig:9}. The one in Fig. \ref{fig:9}a was the first one ever as obtained at ($r_0=1$; $r_1=1$; $r_2=2.26$). Calculated radii are ($r_3=0.81217036$; $r_4=1.75470801$; $r_5=0.39249703$; $r_6=1.36072028$; $r_7=0.10652797$; $r_8=0.22963098$). The other configuration in Fig. \ref{fig:9}b is obtained at ($r_0=1$; $r_1=0.8$;$r_2=2.1$) with calculated radii ($r_3=2.12301926$; $r_4=2.65756162$; $r_5=0.2340634$; $r_6=2.69771176$; $r_7=0.21976258$; $r_8=0.41771213$).  The latter configuration is unique because it has the least ratio of  $12.27$ between the largest and the smallest radii in the radii set of all the $587$ sets calculated. The point ($r_1=0.8$;$ r_2=2.1$) which corresponds to this unique configuration is marked with a star inside the domain of the allowed radii $r_1$ and $r_2$ that we explored and depicted in Fig. \ref{fig:9}c. In principle, the rest $6$ radii can be visualized as $6$ surfaces roofing this domain. The boundary points of the domain may include configurations where the radius of at least one of the cylinder either expands to infinity (mostly the upper part of the domain as shown in Fig. \ref{fig:9}d) or shrinks to zero (mostly the lower part of the domain).

The topological classification of 9-knot is possible in analogy with 7-knot discussed above. For example, exactly like in the 7-knot with equal cylinders only three cylinders: 6th (violet), 3rd (green), and 2nd (blue) are unknotted for 9-knot. The most knotted "inner" cylinders can be released only after removing $2$ others one by one. For the matrix $P_{ij}^{(9)}$ we obtained for the three configurations shown in Fig. \ref{fig:9}

\begin{equation}
P_{ij}^{(9)}=\left( \begin{array}{ccccccccc}
0&+1&+1&+1&+1&+1&+1&+1&+1\\
+1&0&+1&+1&+1&-1&+1&+1&-1\\
+1&+1&0&-1&-1&-1&+1&+1&-1\\
+1&+1&-1&0&-1&+1&+1&-1&-1\\
+1&+1&-1&-1&0&+1&-1&-1&+1\\
+1&-1&-1&+1&+1&0&+1&-1&-1\\
+1&+1&+1&+1&-1&+1&0&-1&-1\\
+1&+1&+1&-1&-1&-1&-1&0&-1\\
+1&-1&-1&-1&+1&-1&-1&-1&0
\end{array} \right)
\label{eq:20}
\end{equation}
with its determinant $|P_{ij}^{(9)} |=0$. This degeneracy may point out to some additional symmetry hidden in this ultimate configuration of 9-knot. We did not check for all $587$ configurations calculated for the outline of the domain in Fig. 9c but where we checked the matrix was the same. Thus we conclude that the topology was likely to be preserved. One may notice that the left upper corner $8\times 8$ block of the matrix exactly coincides with Eq. (\ref{eq:18}) which was the consequence of constructing the 9-knot by introducing the ninth cylinder into the previous 8-knot configuration and then recalculating.  It is for now an open question whether other topologies exist for 9-knot.

\section{Conclusion}
\label{sec:6}
We have given only a glimpse of a rich world of geometry of cylinder configurations that is to be revealed in detail. Numerous generalizations and questions come to mind, such as: what if the cylinders are having arbitrary cross-sections: elliptic, polygonal, or just a segment, a cylinder being degenerated into a flat band (excluding though the degenerate solutions of more than two edges touching in just one point)? What if such bands are not flat, what if the axes of the cylinders are arbitrary curves that still connect "one infinity" with the other one, etc.? With the help of our method of calculating the degrees of freedom it is possible to estimate the relevant numbers but there should be advanced analytics or computations to make sure whether an arrangement is realizable like it happened in case of 7-, 8-, and 9-knots.

In our previous works it was already shown that entanglement of the cylinders brings unconventional mechanical (auxetic) properties to the ensembles of cylinders by using 6-knots as building blocks. Now we introduced a matrix that identifies the topology of cylinder configurations and probably may be used for a description of ensembles of entangled cylinders or rods that goes beyond the usual Onsager approximation of low packing density and pairwise interaction of rigid rods. On this way one may expect more complicated entangled and thermodynamically stable structures than the well-known nematic ordering of rigid rods to exist, for example, similar to spin-glasses.

Many objects in the field theory such as strings and vortexes have the geometry of cylinders. The difference between the cylinder and the sphere may resemble the difference between particles with spin and spinless particles. One may notice an intriguing parallel between the geometry of infinite cylinders, punching the "celestial sphere", and twistor geometry of light rays \cite{Ref5}. A light ray in this picture lies along the axis which is not structureless but is equipped with the information of the helicity (chirality). Contacting cylinders and other cylinder-like geometrical constructions may provide a complementary picture for twistors.  Note that as we showed above just a configuration of two cylinders produces chirality similar to the helicity of the twistor \cite{Ref5}.

We found here that the topology of entanglement of cylinders can be described in easier (and one to one) manner than for conventional one-component knots. This finding may be useful for the Information Theory as well.

We anticipate that some confusion with the terminology may occur in future. There has been a term of n-knot which is a single n-dimensional sphere $S^n$ embedded in m-dimensional sphere $S^m$. Therefore it might be needed to alter our term "n-knot" for a configuration of n mutually touching cylinders for something that we may suggest to be "n-cross".

\begin{acknowledgements}
One of the authors (PVP) thanks Stanislaw Pikhitsa for help in calculating 587 configurations for 9-knots. 
\end{acknowledgements}



\end{document}